\documentclass[12pt]{article}
\usepackage{amssymb}
\usepackage{latexsym}

\begin{document}

\newtheorem{theorem}{Theorem}
\newtheorem{lemma}{Lemma}
\newtheorem{proposition}{Proposition}
\newtheorem{corollary}{Corollary}
\newtheorem{remark}{Remark}

\title{Reduction theory for mapping class groups \\
and applications to moduli spaces}
\author{Enrico Leuzinger}
\date{January 10, 2008}
\maketitle

\begin{abstract} 
Let  $S=S_{g,p}$  be a compact, orientable surface of genus $g$ with $p$ punctures and
such that  $d(S):=3g-3+p>0$. The mapping class group $\textup{Mod}_S$ acts properly discontinuously on the Teichm\"uller space $\mathcal T(S)$ of marked hyperbolic structures on $S$.
The resulting quotient $\mathcal M(S)$ is the moduli space of isometry classes of hyperbolic  surfaces.  We provide a version of precise reduction theory for finite index subgroups of $\textup{Mod}_S$, i.e.,
a description of exact fundamental domains.
As an application we show that the asymptotic cone of the moduli space $\mathcal M(S)$ 
endowed with the Teichm\"uller metric  is bi-Lipschitz equivalent  to
the Euclidean cone over the finite simplicial (orbi-) complex $ \textup{Mod}_S\backslash\mathcal C(S)$, where $\mathcal C(S)$ of $S$ 
is the complex of curves of $S$. 
We also show that  if $d(S)\geq 2$, then $\mathcal M(S)$ does \emph{not} admit a finite volume Riemannian metric of (uniformly bounded) positive scalar 
curvature  in the bi-Lipschitz class of the Teichm\"uller metric. These two applications  confirm conjectures of Farb.
\end{abstract}

\vspace{6ex}



\maketitle

\section{Introduction and main results}

The chief goal of this article is to study the large-scale geometry of moduli spaces of Riemann
surfaces endowed with the Teichm\"uller metric. In this introductory section we recall basic notions and present an outline of our main results.

Let $S=S_{g,p}$ be a compact, orientable surface  of genus $g$ with $p$ punctures
such that $3g-3+p>0$. This last assumption implies that $S$ carries  finite volume Riemannian metrics of constant  curvature $-1$ and $p$  cusps. 
A {\it marked hyperbolic surface} is a pair $(X, [f])$ where $X$ is a smooth surface equipped with a complete Riemannian metric of constant curvature $-1$ and where $[f]$ denotes the isotopy class of  a diffeomorphism   $f:X\longrightarrow S$
mapping cusps to punctures.
Two marked surfaces 
$(X_1,[f_1])$ and
$(X_2,[f_2])$ are equivalent
if there is an isometry  $h:X_1\longrightarrow X_2$ such that
$[f_2\circ h]=[f_1]$.  The collection of these equivalence classes is (one possible definition of) \emph{the Teichm\"uller space}
$\mathcal{T}(S)$ of $S$. The corresponding {\it moduli space} $\mathcal M(S)$ of isometry classes of hyperbolic surfaces is obtained by forgetting the marking. More precisely, consider the
 {\it mapping class group} $\textup{Mod}_S$, i.e., the group of all orientation preserving 
 diffeomorphisms of $S$, which fix the punctures,
modulo isotopies which also fix the punctures. Then $\textup{Mod}_S$ acts on $\mathcal T(S)$ 
according to the formula $h\cdot x=h\cdot (X,[f]):=(X,[h\circ f])$, for $ h\in \textup{Mod}_S, x\in \mathcal T(S)$, and  one has $\mathcal M(S)=\Gamma\backslash \mathcal T(S)$.

The 
\emph{complex of curves} $\mathcal C(S)$ of $S$ is an infinite (even locally infinite) simplicial complex
of dimension $d(S)-1$.
It was introduced by  Harvey in \cite{Har1} as an analogon in the context of Teichm\"uller spaces of the (rational) Tits building associated
to an  arithmetic group.
The \emph{vertices} of $\mathcal C(S)$  are the isotopy classes of simple closed curves (called \emph{circles}) on $S$, which are non-trivial (i.e., not contractible in $S$ to a point or to a component of $\partial S$). We denote the isotopy class
of a circle $C$ by $\langle C\rangle$. A set of $k+1$ vertices $\{\alpha_0,\ldots \alpha_k\}$ span a
$k$-\emph{simplex} of $\mathcal C(S)$ if and only if $\alpha_0=\langle C_0\rangle,\ldots,\alpha_k=\langle C_k\rangle$
for some set of pairwise non-intersecting circles $C_0,\ldots, C_k$. 
For a simplex $\sigma\in \mathcal C(S)$ we denote by $|\sigma|$ the number of its vertices. A crucial fact  is that $\mathcal C(S)$ is a thick chamber complex, i.e., every simplex is the face of a maximal simplex.  Moreover the mapping class group $\textup{Mod}_S$ acts simplicially on $\mathcal C(S)$ and the quotient  $ \textup{Mod}_S\backslash\mathcal C(S)$ is a finite (orbi-) complex  (see \cite{Har1}, Proposition 1).

In order to simplify  exposition and proofs of the present  article we work  in the framework of manifolds and simplicial complexes rather than orbifolds and orbicomplexes. We thus consider finite index, torsion-free subgroups  $\Gamma$ of $ \textup{Mod}_S$ which in addition consist of pure mapping classes.
Recall that a maping class $h\in \textup{Mod}_S$ is called {\it pure} if it can be represented by a diffeomorphism $f:S\longrightarrow S$ fixing (pointwise) some union $\Lambda$ of disjoint and pairwise non-isotopic nontrivial circles on $S$ and such that $f$ does not permute the components of $S\setminus \Lambda$ and induces on each component of the cut surface $S_{\Lambda}$ a diffeomorphism isotopic to a pseudo-Anosov
or to the identity diffeomorphism (see \cite{Iv2}, \S\,7.1). It is well-known that such subgroups $\Gamma$  exist. For example one can take $\Gamma=\Gamma_S(m)$,  the kernel of the natural homomorphism
$\textup{Mod}_S\longrightarrow \textup{Aut}(H_1(S,\mathbb Z/m\mathbb Z)), m\geq 3$, defined by the action of diffeomorphisms on homology (see e.g. \cite{Iv2}, \S \, 7.1 or \cite{Ha}, \S \, 1.3).

The first step in our aproach to understand the coarse geometry of moduli space 
is the construction of a tiling of Teichm\"uller space $\mathcal T(S)$, i.e., a $\Gamma$-invariant decomposition into disjoint subsets, which in turn yields a tiling of  moduli space $\mathcal M(S)$.
In order to formulate that result, we need  {\it length functions}.
Let $\alpha$ be a vertex of $\mathcal C(S)$, i.e., $\alpha=\langle C\rangle$
for a (non-trivial) circle $C$. 
Since $d(S)>0$,  any point $x\in \mathcal T(S)$ represents a finite volume Riemann surface of  curvature $-1$ with $p$ cusps. On the surface $x$ the isotopy class of $\alpha$ 
contains a unique geodesic loop; let $l_{\alpha}(x)$ denote its length. This defines 
a (smooth) function $l_{\alpha}: \mathcal T(S)\longrightarrow \mathbb R_{>0}$ for every vertex $\alpha\in \mathcal C(S)$.

For $\varepsilon>0$ we then define the {\it $\varepsilon$-thick part of Teichm\"uller space} 
$$\textup{Thick}_{\varepsilon}\mathcal T(S):=
\{x\in\mathcal T(S)\mid l_{\alpha}(x)\geq \varepsilon, \forall \alpha\in\mathcal C(S)\}.
$$
 This set is 
$\Gamma$-invariant and its quotient is the {\it $\varepsilon$-thick part $\textup{Thick}_{\varepsilon}\mathcal M(S)$ of moduli space}. 

We can now state  a concise version of 
our first main result.

\vspace{2ex}

\noindent{\bf Theorem A} \ {\bf(Tiling of moduli space)} {\it Let $S=S_{g,p}$ be a compact, orientable surface of genus $g$ with $p$ punctures such that $3g-3+p>0$. Let $\Gamma$ be   a torsion free, finite index subgroup of the mapping class group of $S$ consisting of pure elements. Further let  $\mathcal M(S)= 
\Gamma\backslash\mathcal T(S)$ be the corresponding moduli space of 
Riemann surfaces and let $\mathcal E$ be the set of simplices of the  finite complex $ \Gamma\backslash\mathcal C(S)$. Then there is $\varepsilon =\varepsilon (S)>0$ such that the $\varepsilon$-thick part  $\textup{Thick}_{\varepsilon}\mathcal M(S)$ of moduli space is a compact submanifold with corners whose boundary is the union of compact sets 
$ \mathcal B_{\varepsilon}(\sigma)$ indexed by $\mathcal E$ and
 this structure extends to a thick-thin tiling of the entire moduli space: there is a \textup{disjoint union}
$$
\mathcal M(S)= \textup{Thick}_{\varepsilon}\mathcal M(S)\sqcup\bigsqcup_{\sigma\in \mathcal E}\textup{Thin}_{\varepsilon}(\mathcal M(S),\sigma)
$$ 
where the $(\varepsilon,\sigma)$-thin part $\textup{Thin}_{\varepsilon}(\mathcal M(S),\sigma)$ is diffeomorphic to 
 $\mathcal B_{\varepsilon}(\sigma)\times \mathbb R_{>0}^ {|\sigma|}$.}

\vspace{2ex}

For more details on the structure of the boundary faces  
$ \mathcal B_{\varepsilon}(\sigma)$ see Proposition \ref{boundary} in Section 2.1 below.
Theorem A (which we prove in Section 2) can be considered as a version of precise reduction theory for mapping class groups.
In fact, there are striking parallels to  precise reduction theory for arithmetic groups as developped e.g. by 
Langlands \cite{La}, Osborne-Warner \cite{OW} and Saper \cite{Sa} (see also \cite{BJ}).
In Section 2 we further emphasize this analogy  using the notion of   {\it parabolic subgroups} of $\Gamma\subset\textup{Mod}_S$ (compare Corollary \ref{para}).

In Section 3 we study (coarse) metric properties of the tiling described in Theorem A. 
We in particular  construct a quasi-isometric  approximation  of  moduli
space $\mathcal M(S)=\textup{Mod}_S\backslash\mathcal T(S)$ endowed with the distance function $d_{\mathcal M}$ 
induced by the  Teichm\"uller metric. This approximation is provided by the  {\it
asymptotic cone (or tangent cone at infinity)} defined as the Gromov-Hausdorff limit of rescaled pointed
metric spaces:
$$ \textup{As-Cone}(\mathcal M(S)) :=  {\mathcal H}-{\lim}_{n\rightarrow\infty} (\mathcal M(S), x_0, \frac{1}{n} d_{\mathcal M}),$$ 
where $x_0$ is an arbitrarily chosen  point of $\mathcal M(S)$.
  We remark that
in contrast to the case considered here, the   definition of an asymptotic
 cone in general involves the use of ultrafilters, and the limit space may   depend on the
chosen
ultrafilter. Various aspects of  asymptotic cones of general spaces are discussed in  Gromov's essay \cite{Gr}.
In some cases asymptotic cones are  easy to describe. For
example,
if  $V$ is a finite volume Riemannian manifold of   strictly negative sectional curvature and with $k$ cusps,  then   $\textup{As-Cone}(V)$ is
a ``cone''  over $k$ points, i.e., $k$ rays with a common origin.  For  a Riemannian product $V= V_1\times V_2$, where
$V_1$, $V_2$ are as in the previous example and each has only one cusp,
$\textup{As-Cone}(V)$  can be identified with the first quadrant in $\mathbb R^2$.
Much more intricate are quotients of $SL(n,\mathbb R)/SO(n)$ by congruence subgroups of
$SL_n(\mathbb Z)$ or more general locally symmetric spaces $V=\Gamma\backslash G/K$ of higher rank.
For such $V$  $\textup{As-Cone}(V)$ is \emph{isometric} to the Euclidean cone over the finite simplicial complex 
 given by the quotient of the \emph{rational} Tits building of $G$ modulo $\Gamma$
(see \cite{Hat}, \cite{L2}, \cite{L3}).

\vspace{1ex}

 Theorem A allows us to determine the asymptotic cone of moduli space. Notice the striking similarity of  Theorem B below with the result for locally symmetric 
spaces  described above.

\vspace{2ex}

\noindent{\bf Theorem B} \textup{\bf{(Asymptotic cone of moduli space)}}\  {\it Let $S$ and $\Gamma$ be   as in Theorem A.  Further let the moduli space $\mathcal M(S)= 
\Gamma\backslash\mathcal T(S)$  of 
Riemann surfaces be endowed with the metric $d_{\mathcal M}$ induced by the Teichm\"uller metric. 
 Then the
 asymptotic cone of $(\mathcal M(S), d_{\mathcal M})$  is bi-Lipschitz equivalent  to
the Euclidean cone over the finite simplicial complex $\Gamma\backslash|\mathcal C(S)|$, where 
$|\mathcal C(S)|$ is (a geometric realization of) the complex of curves of $S$. 
In particular, $\dim \textup{As-Cone}(\mathcal M(S)) =\dim |\mathcal C(S)|+1=3g-3+p=\frac{1}{2}\dim \mathcal T(S)$.}

\vspace{2ex}

A variant of Theorem B also holds for the orbifold moduli space $\textup{Mod}_S\backslash\mathcal T(S)$, which is finitely covered by
$\Gamma\backslash\mathcal T(S)$.
This confirms (the strong version of) a conjecture of Farb (\cite{Fa}, Conjecture 4.7). 
The proof of Theorem B, which we give in Section 3, is based on two key ingredients. The first one is McMullen's modification of 
the Weil-Petersson metric \cite{Mc}. In fact,
 we only need information on the bi-Lipschitz class of the Teichm\"uller metric and McMullen's
metric belongs to that class. We remark in passing that there are several other complete metrics
comparable to the Teichm\"uller metric (see \cite{LSY}). The second ingredient of the proof of Theorem B  is an expansion of the WP-metric on $(\varepsilon,\sigma)$-thin parts $\textup{Thin}_{\varepsilon}(\sigma;\mathcal M(S))$ due to Wolpert \cite{Wo2}.

\vspace{1ex}

\noindent {\bf Remark}. For an
arithmetic lattice $\Gamma$ in   a semisimple Lie group $G$ with associated locally symmetric space $\Gamma\backslash G/K$ one has  
$\dim \, \textup{As-Cone} (\Gamma\backslash G/K)=\mathbb Q$-$\textup{rank}\, \Gamma$ (see \cite{L1}). On the other hand,
 by Theorem A,  $\dim \, \textup{As-Cone} (\mathcal M(S))=3g -3+p$. The number $d(S):= 3g-3+p$,
which measures the topological complexity of $S$, might  thus 
 be considered as ``the'' $\mathbb Q$-$\textup{rank}$ of the mapping class group;
 compare the discussion in \cite{Fa}, \S  \ 4.2. The \emph{geometric rank } of $\mathcal T(S)$
is the maximal dimension of a quasi-isometrically embedded  Euclidean space.  In \cite{BM} it is shown
  that the  \emph{Weil-Petersson} geometric rank of $\mathcal T(S)$ equals 
  $[\frac{d(S)+1}{2}]$. Thus  $\textup{WP-rank} \, \mathcal T(S)\le d(S)= \mathbb Q$-$\textup{rank}\, \textup{Mod}_S$.
  This is in contrast to the case of locally symmetric spaces where $\textup{rank} \, G/K= \mathbb R$-$\textup{rank}\, G\geq  \mathbb Q$-$\textup{rank}\, \Gamma$ and rises the question if the {\it Teichm\"uller} geometric rank 
dominates $\mathbb Q$-$\textup{rank}\, \textup{Mod}_S$. Notice that, by Theorem B, there are Euclidean balls of arbitrarily large
radius quasi-isometrically embedded in $\mathcal T(S)$ and $\mathcal M(S)$.

\vspace{1ex}
 
 Theorem A together with subresults of its proof also provides an obstruction to a Riemannian metric of  positive scalar curvature
 in the quasi-isometry class of the Teichm\"uller metric.

\vspace{2ex}

\noindent{\bf Theorem C}\ \textup{\bf {(No positive scalar curvature)}} \  {\it Let $S=S_{g,p}$ be a compact, orientable surface with genus $g$ and $p$ punctures such that $d(S)= 3g-3+p\geq 2$ (i.e., $\mathbb Q$-$\textup{rank} \ \textup{Mod}_S\geq 2$). Further let $\Gamma\subset \textup{Mod}_S$ be   a torsion free, finite index subgroup. 
Then the corresponding moduli space $\mathcal M(S)= 
\Gamma\backslash\mathcal T(S)$ does not admit a finite volume Riemannian metric 
of uniformly bounded
positive scalar curvature in the quasi-isometry class of the Teichm\"uller metric.}

\vspace{2ex}

This result has also been conjectured by Farb (\cite{Fa}, Conjecture 4.6).
For locally symmetric spaces associated to arithmetic groups of $\mathbb Q$-$\textup{rank} \geq 2$ the analogon of Theorem C was
proven by  Chang \cite{Ch}. We show in Section 4 how the ideas and phenomena underlying his proof can be extended  to the present setting of Teichm\"uller and moduli spaces.

\vspace{1ex}

\noindent{\bf Notation}. We write $A\asymp B$ (resp. $A\succ B$) if there is a constant $C>0$
such that $C^{-1}A<B<CA$ (resp. $CA>B$).

\section{Reduction theory for mapping class groups}

The goal of this section is to prove Theorem A of Section 1.

\subsection{The complex of curves and  submanifolds with corners}

The following two lemmata will  be used frequently.

\begin{lemma}\label{length} Let $\alpha\in\mathcal C(S)$ and $x\in \mathcal T(S)$. Then for any $h\in \textup{Mod}_S$ holds
$$
l_{\alpha}(h\cdot x)=l_{h^{-1}\cdot \alpha}(x).
$$
\end{lemma}

{\it Proof}. If $x=(R,[f])$, with $f:R\longrightarrow S$,  then $h\cdot x=(R,[h\circ f])$. By definition,  $l_{\alpha}(x)$ is the hyperbolic length of the unique closed geodesic
in the isotopy class of $f^{-1}(\alpha)$ in $R$: $l_{\alpha}(x)=L_{hyp}([f^{-1}(\alpha)])$. Hence
$
l_{\alpha}(h\cdot x)=L_{hyp}([f^{-1}\circ h^{-1}(\alpha)])=l_{h^{-1}\cdot \alpha}(x).
$
\hfill $\Box$

\begin{lemma} \label{disjoint}There exists a universal constant $c=c(S)$ such that the following holds.
If $\alpha_i$,  $0\leq i\leq k$,
are simple geodesic loops on a hyperblic surface $x\in \mathcal T(S)$ such that $l_{\alpha_i}(x)<c$ for all $0\leq i\leq k$, then the loops $\alpha_i$ are disjoint and hence define a $k$-simplex in $\mathcal C(S)$.

\end{lemma}

For a proof of Lemma \ref{disjoint} see \cite{Abi}, \S\,  3.3. 

\vspace{2ex}

Recall that for $\varepsilon >0$ the {\it $\varepsilon$-thick part of Teichm\"uller space} is the set 
$\textup{Thick}_{\varepsilon}\mathcal T(S):=\bigcap \{x\in\mathcal T(S)\mid l_{\alpha}(x)\geq \varepsilon \ \textup{for all vertices}\  \alpha\in \mathcal C(S)\}$. The family of closed 
sets $\{x\in\mathcal T(S)\mid l_{\alpha}(x)\leq \varepsilon\}$ turns out to be locally finite if $\varepsilon$ is sufficiently small. 
 For every vertex $\alpha\in \mathcal C(S)$ we further set $\mathcal H_{\varepsilon}(\alpha):=\{x\in \textup{Thick}_{\varepsilon}\mathcal T(S)\mid l_{\alpha}(x)=\varepsilon\}$.
Then (for fixed $\varepsilon$ sufficiently small) we have $\partial \textup{Thick}_{\varepsilon}\mathcal T(S)=\bigcup_{\alpha} \mathcal H_{\varepsilon}(\alpha)$.

The following proposition is due to Ivanov (see \cite{Iv1}, \S\, 4.6 and also \cite{Ha}, Ch. 3).

\begin{proposition}\label{corners} There is $\varepsilon_0 >0$ depending only on $S$ such that for
all $\varepsilon\leq \varepsilon_0$ the following holds:
\textup{(1)} \ The $\varepsilon$-thick part  $\textup{Thick}_{\varepsilon}\mathcal T(S)$ is a submanifold of $\mathcal T(S)$ with corners,
i.e., locally modelled on a cube $[0,1]^n\subset \mathbb R^n$, and invariant under $\Gamma$.

\textup{(2)} \ The  complex of curves\  $\mathcal C(S)$ is the nerve of the
covering of $\partial \textup{Thick}_{\varepsilon}\mathcal T(S)$ by the closed sets $\mathcal H_{\varepsilon}(\alpha)$.
In particular there is a one-to-one correspondence between the  boundary faces  $\bigcap_{\alpha\in\sigma}\mathcal H_{\varepsilon}(\alpha)$ of $\partial \textup{Thick}_{\varepsilon}\mathcal T(S)$ and the simplices of $\mathcal C(S)$.

\textup{(3)} \  For  $\sigma\in\mathcal C(S)$ the  boundary face $\bigcap_{\alpha\in\sigma}\mathcal H_{\varepsilon}(\alpha)$
in $\partial \textup{Thick}_{\varepsilon}\mathcal T(S)$ is diffeomorphic to $\textup{Thick}_{\varepsilon}\mathcal T(S_{\sigma})\times \mathbb R^{|\sigma|}$, where $S_{\sigma}$ is the 
result of cutting $S$ along (non-intersecting) circles from the isotopy classes  $\alpha\in\sigma$ and such
that the length  of each boundary circle is $\varepsilon$.
\end{proposition}

We emphasize the following simple consequence of Proposition \ref{corners}\ (2).

\begin{lemma} \label{bdpoint-simplex} To any point $x\in \partial \textup{Thick}_{\varepsilon}\mathcal T(S)$ there 
corresponds a unique simplex $\sigma\in \mathcal C(S)$ such that $x\in \bigcap_{\alpha\in\sigma}\mathcal H_{\varepsilon}(\alpha)$ and $\sigma$ is maximal for this condition.
\end{lemma}

{\it Proof}. In view of Proposition \ref{corners}\ (2) it suffices to observe that $\bigcap_{\alpha\in\sigma}\mathcal H_{\varepsilon}(\alpha)\supset \bigcap_{\alpha\in\tau}\mathcal H_{\varepsilon}(\alpha)$
if and only if $\sigma\subset \tau$.
\hfill$\Box$

\vspace{2ex}

Recall from Section 1 that instead of the mapping class group $\textup{Mod}_S$ itself we consider 
 finite index subgroups $\Gamma$, which are torsion free and in addition consist of pure mapping classes.
This allows us to eliminate all difficulties related to elements of finite order and to possible
permutations of components of cut surfaces by reducible diffeomorphisms.

A subgroup $\Gamma\subset \textup{Mod}_S$ with the above properties leaves invariant 
the $\varepsilon$-thick part  $\textup{Thick}_{\varepsilon}\mathcal T(S)$ of Teichm\"uller space
 (by Lemma \ref{length}) and acts freely on it.
The corresponding   quotient
$\textup{Thick}_{\varepsilon}\mathcal M(S):=\Gamma\backslash \textup{Thick}_{\varepsilon}\mathcal T(S)$, which we call the \emph{$\varepsilon$-thick part of moduli space}, is compact (see \cite{Iv1}, \S\, 4). This yields  the following corollary.  For an analogous result  for locally symmetric spaces
see \cite{L1}.

\begin{corollary}  \label{exhaustion} There exists  $\varepsilon_0>0$ depending only on $S$, such that there is a $\Gamma$-invariant exhaustion of Teichm\"uller space $\mathcal T (S) = \bigcup_{\varepsilon\leq \varepsilon_0}\textup{Thick}_{\varepsilon}\mathcal T(S)$, which induces an exhaustion of moduli space 
$\mathcal M(S)$  by  polyhedra, i.e., by compact submanifolds  with corners: $\mathcal M (S) = \bigcup_{\varepsilon\leq \varepsilon_0}\textup{Thick}_{\varepsilon}\mathcal M(S)$.
\end{corollary}

We choose $\varepsilon\leq \varepsilon_0$, where $\varepsilon_0$ is as in Proposition \ref{corners}.
By Proposition \ref{corners}(2) there is a one-to-one correspondence between the boundary faces of $\partial \textup{Thick}_{\varepsilon}\mathcal T(S)$ and the simplices of the complex of curves $\mathcal C(S)$. Clearly, the torsion-free finite index subgroup $\Gamma$ of the mapping class group  acts equivariantly with respect to that correspondence. 
The quotient 
$\Gamma\backslash \mathcal C(S)$ is a finite simplicial complex  (see \cite{Har1}, Proposition 1). We denote by $\mathcal E$ 
the  set of simplices of  
$\Gamma\backslash \mathcal C(S)$ or, equivalently, the set of the boundary faces of $\partial\textup{Thick}_{\varepsilon}\mathcal M(S)=\Gamma\backslash \partial \textup{Thick}_{\varepsilon}\mathcal T(S)$. The next proposition describes the fine structure of the boundary $\partial\textup{Thick}_{\varepsilon}\mathcal M(S)$.

\begin{proposition} \label{boundary}
Let $\varepsilon_0>0$  be as in Proposition \ref{exhaustion}.  Then, for each $\varepsilon\leq \varepsilon_0$, the boundary $\partial \textup{Thick}_{\varepsilon}\mathcal M(S)$ consists of  a finite
number of
 faces  $\mathcal B_{\varepsilon}(\sigma),\ \sigma\in {\mathcal E}$.  Each $\mathcal B_{\varepsilon}(\sigma)$ is
 a (trivial) torus bundle over (a finite covering of) the $\varepsilon$-thick part of the
moduli space of the cut surface $S_{\sigma}$ (i.e., the surface obtained by cutting $S$ along the circles of $\sigma$):
$$
0\rightarrow T^{|\sigma|}\rightarrow \mathcal B_{\varepsilon}(\sigma)
\rightarrow \textup{Thick}_{\varepsilon}\mathcal M(S_{\sigma})\rightarrow 0.
$$
Moreover, the nerve of the covering of $\partial \textup{Thick}_{\varepsilon}\mathcal M(S)$ by these faces is isomorphic to the finite simplicial complex
$\Gamma\backslash \mathcal C(S)$. 
For  a  simplex $\tau$ of maximal dimension $d(S)$, $\mathcal M(S_{\tau})$ and hence  $\textup{Thick}_{\varepsilon}\mathcal M(S_{\tau})$
 reduces to a point and $\mathcal B_{\varepsilon}(\tau)$ is a $d(S)$-dimensional torus.

\end{proposition}

{\it Proof}.  Given Proposition \ref{corners}, it remains to determine the stabilizer in $\Gamma$ of a
boundary face $\bigcap_{\alpha\in \sigma} \mathcal H_{\varepsilon}(\alpha)$ of $\partial \textup{Thick}_{\varepsilon}\mathcal T(S)$. Thus, we assume that for $x$ in  $\bigcap_{\alpha\in \sigma} \mathcal H_{\varepsilon}(\alpha)$ (with $\sigma$ maximal as in Lemma \ref{bdpoint-simplex}) there is $h\in \Gamma$ such that $h\cdot x$  also is in 
$\bigcap_{\alpha\in \sigma} \mathcal H_{\varepsilon}(\alpha)$. Since $l_{h^{-1}\cdot \alpha}(x)=l_{\alpha}(h\cdot x)$ for all $\alpha\in \sigma$ by Lemma \ref{length}, we have $x\in \bigcap_{\beta\in h^{-1}\cdot\sigma\vee\sigma} \mathcal H_{\varepsilon}(\beta)$ and as $\sigma$ is maximal for this condition $h$ must stabilize $\sigma$: $h\cdot\sigma=\sigma$. 

Now by \cite{Har2}, \S \, 4,  in general the stabilizer of $\sigma$ in $\textup{Mod}_S$ is a group extension of a finite
group (a  subgroup of the automorphism group, say $A(\sigma)$, of the decomposition graph associated to the set of circles corresponding  to $\sigma$) over  the normalizer of the twist group $\textup{Tw}(\sigma)$ of $\sigma$, i.e., the free abelian
group  generated by the $|\sigma|$ Dehn twists $\tau_{\alpha}$ of circles of $\sigma$. Since $\Gamma$ consists of pure mapping classes we have $\Gamma\cap A(\sigma)=\{\textup{id}\}$. Moreover,  we also have $h^{-1}\cdot \alpha=\alpha$ for all $\alpha\in \sigma$ and thus $h^{-1}\tau_{\alpha}h=\tau_{h^{-1}\cdot \alpha}=\tau_{\alpha}$ for  all $\alpha\in \sigma$. Hence $h\in \textup{stab}_{\Gamma}(\sigma)=Z_{\Gamma}(\textup{Tw}(\sigma))$.

Finally, in order to determine $Z_{\Gamma}(\textup{Tw}(\sigma))$,  note that a  pure mapping class $h$ which commutes with a Dehn twist ``of $\sigma$" preserves setwise every circle of $\sigma$ and every part
of   $S\setminus \sigma$ (see \cite{Iv2}, 7.5).  Thus $h$ belongs to the direct product $(\textup{Mod}_{S_{\sigma}}\cap\Gamma)\times \textup{Tw}(\sigma)$ and the result follows.
\hfill $\Box$

\subsection{Tilings of Teichm\"uller  spaces}

Recall that the complex of curves $\mathcal C(S)$ is a thick chamber complex.
Given a simplex $\sigma\in \mathcal C(S)$ we can thus choose a simplex $\tau\in \mathcal C(S)$ of maximal dimension 
$d(S)-1=3g-4+p$ containing $\sigma$. Then there are  \emph{adapted Fenchel-Nielsen coordinates} 
on Teichm\"uller space $\mathcal T(S)$  (see e.g. \cite{Abi}), i.e., a diffeomorphism
$$
\Phi_{\tau}: \mathcal T(S)\longrightarrow \mathcal T(S_{\sigma})\times \mathbb R^{|\sigma|}\times \mathbb R_{ > 0}^{|\sigma|}\ ;\ p\longmapsto (s(p), \theta(p),l(p)),
$$
where $s=(\theta_{\beta},l_{\beta})_{\beta\in\tau\setminus\sigma}$ parametrizes $\mathcal T(S_{\sigma})$, 
$\theta=(\theta_{\alpha})_{\alpha\in\sigma}$ are twist parameters on $\mathbb R^{|\sigma|}$ and
$l= (l_{\alpha})_{\alpha\in\sigma}$ are coordinates on $\mathbb R_{>
 0}^{|\sigma|}$ 
(here and elsewhere $|\sigma|$ denotes the number of vertices of $\sigma$).

Consider $x\in \partial \textup{Thick}_{\varepsilon}\mathcal T(S)$. By Lemma \ref{bdpoint-simplex}, $x$ determines a unique (largest) simplex $\sigma\in\mathcal C(S)$. Let  $\tau\in \mathcal C(S)$ be a simplex of maximal dimension containing $\sigma$ and let  $\Phi_{\tau}$ be  adapted FN-coordinates.
We define the \emph{outer cone } at $x\in\partial \textup{Thick}_{\varepsilon}\mathcal T(S)$  as the preimage 
$$CO(x):=\Phi_{\tau}^{-1}\{(s(x),  \theta(x),(l_{\alpha})_{\alpha\in \sigma})\mid l_{\alpha}    <     \varepsilon\ \ \textup{for all}\ \      \alpha\in \sigma\}.
$$

Note that $CO(x)$ is diffeomorphic to the open {\it hyperoctant} $\mathbb R_{>0}^{|\sigma|}$.

\begin{lemma} \label{outercone} 
This definition of an outer cone $CO(x)$  is independent of the chosen simplex $\tau$ of maximal dimension (resp. its adapted FN-coordinates $\Phi_{\tau}$). 
Furthermore, the set of all outer cones
of $\{CO(x)\mid x\in \partial \textup{Thick}_{\varepsilon}\mathcal T(S)\}$ is $\Gamma$-invariant.
\end{lemma}

{\it Proof}.\  If $\tilde{\tau}$ is another  simplex 
of maximal dimension in $\mathcal C(S)$ that also contains $\sigma$ we have a diffeomorphism (adapted FN-coordinates)
$$\Phi_{\tilde{\tau}}: \mathcal T(S)\longrightarrow \mathcal T(S_{\sigma})\times \mathbb R^{|\sigma|}\times \mathbb R_{ > 0} ^{|\sigma|}\ ;\ p\longmapsto (\tilde{s}(p),  \tilde{\theta}(p), (l_{\alpha}(p))_{\alpha\in \sigma}).
$$
Hence
$$\Phi_{\tau}\circ \Phi_{\tilde{\tau}}^{-1}((\tilde{s}(p),   \tilde{\theta}(p),(l_{\alpha}(p))_{\alpha\in \sigma})=
(s(p), \theta(p), (l_{\alpha}(p))_{\alpha\in \sigma}), \ \ p\in \mathcal T(S).
$$ 
Since $(l_{\alpha})_{\alpha\in \sigma}$ are coordinate functions (for both $\Phi_{\tilde{\tau}}$ {\it and}
$\Phi_{\tau}$) the implicit function theorem implies that there is a diffeomorphism
$$
\varphi: \mathcal T(S_{\sigma})\times \mathbb R ^{|\sigma|}\longrightarrow \mathcal T(S_{\sigma})\times \mathbb R ^{|\sigma|};\ \ (\tilde{s}(p),  \tilde{\theta}(p))\longmapsto 
(s(p),   \theta(p)).
$$
In particular we have $(\tilde{s},\tilde{\theta})=\textup{const.}\Leftrightarrow (s,\theta)=\varphi((\tilde{s},\tilde{\theta}))=\textup{const.}$
 This proves the first claim.
The invariance under $\Gamma$ follows from Lemma \ref{length}. \hfill$\Box$

\vspace{2ex}

Any of the   submanifolds with corners $\textup{Thick}_{\varepsilon}\mathcal T(S)$ of Teichm\"uller space   as described 
in Proposition \ref{corners} together with its outer cones yields a tiling (or dissection) of $\mathcal T(S)$ into {\it disjoint} subsets.

\begin{proposition} \label{tileT}  Let $\varepsilon_0 >0$ be as in Proposition \ref{corners} and $\varepsilon\leq \varepsilon_0$.
Then Teichm\"uller space can be written as a $\Gamma$-invariant disjoint union of the $\varepsilon$-thick part and 
all outer cones:
$$
\mathcal T(S)= \textup{Thick}_{\varepsilon}\mathcal T(S)\sqcup \bigsqcup_{x\in\partial\textup{Thick}_{\varepsilon}\mathcal T(S)}CO(x)
= \textup{Thick}_{\varepsilon}\mathcal T(S)\sqcup \bigsqcup_{\sigma\in\mathcal C(S)}\  \bigsqcup_{x\in \bigcap_{\alpha\in\sigma}\mathcal H_{\varepsilon}(\alpha)}CO(x).
$$
\end{proposition}

{\it Proof}.\ Pick $y\in 
\mathcal T(S)\setminus \textup{Thick}_{\varepsilon}\mathcal T(S)$. By Lemma \ref{disjoint} there is a maximal simplex 
$\sigma\in\mathcal C(S)$ (not necessarily of maximal dimension) such that $l_{\alpha}(y)<\varepsilon$ for all $\alpha\in \sigma$.
Choose a simplex $\tau$ of maximal dimension which contains $\sigma$ as a subsimplex.
In adapted FN-coordinates we have $\Phi_{\tau}(y)=(s(y),   \theta(y),(l_{\alpha}(y))_{\alpha\in \sigma})$. By definition of $\sigma$ we have $l_{\beta}(y)\geq\varepsilon$ for all $\beta\in\mathcal C(S)\setminus\sigma$. Thus we find that $x:=(s(y),  \theta(y),(\varepsilon)_{\alpha\in \sigma})\in \partial\textup{Thick}_{\varepsilon}\mathcal T(S)$ and hence
that $y\in CO(x)$ by definition of the latter.

We next show that the outer cones are disjoint. Assume that $z\in CO(x)\cap CO(y)$, $x,y\in \partial\textup{Thick}_{\varepsilon}\mathcal T(S)$.
By Lemma \ref{bdpoint-simplex} there are unique simplices $\sigma$ resp. $\mu$ in $\mathcal C(S)$ associated to $x$ resp. $y$ and,
 by definition of outer cones,  $l_{\alpha}(z)<\varepsilon$ for all $\alpha\in \sigma\cup \mu$.
  By Lemma \ref{disjoint},
  $\sigma\vee\mu$ is a simplex in $\mathcal C(S)$  and thus contained in a simplex of maximal dimension,
say $\tau$. We can thus write
$CO(x)=\Phi_{\tau}^{-1}\{(s(x),   \theta(x),(l_{\alpha})_{\alpha\in \sigma})\mid l_{\alpha}    <     \varepsilon\ \ \textup{for all}\ \      \alpha\in \sigma\}$. In particular we have $l_{\beta}(v)\geq \varepsilon$ for all $v\in CO(x)$ and all  $\beta\in \tau\setminus \sigma$. Our assumption on $z$ thus
implies that $\mu \subseteq \sigma$. Using $CO(y)$ instead of $CO(x)$ we get $\sigma \subseteq \mu$ and hence $\sigma=\mu$.
The formulae for outer cones further yield
$$
(s(x),  \theta(x), (l_{\alpha}(z))_{\alpha\in \sigma})=\Phi_{\tau}(z)=
(s(y),   \theta(y),(l_{\alpha}(z))_{\alpha\in \sigma}).
$$
This implies that  $x=y$ and hence that $CO(x)=CO(y)$ which proves the claim.

In conclusion we have a disjoint union $
\mathcal T(S)= \textup{Thick}_{\varepsilon}\mathcal T(S)\sqcup \bigsqcup_{x\in\partial\textup{Thick}_{\varepsilon}\mathcal T(S)}CO(x)$. The remaining assertion 
then follows from Proposition \ref{corners}.
\hfill$\Box$

\vspace{2ex}

 Given $0<\varepsilon\leq \varepsilon_0$ as in Proposition \ref{tileT} and $\sigma\in\mathcal C(S)$ we define the  
  {\it $(\varepsilon,\sigma)$-thin part} of Teichm\"uller space $\mathcal T(S)$ as a disjoint union of outer cones
$$
\textup{Thin}_{\varepsilon}(\mathcal T(S),\sigma):=
\bigsqcup_{x\in\bigcap_{\alpha\in\sigma}\mathcal H_{\varepsilon}(\alpha)} CO(x).
$$

We can then rewrite the tiling in Proposition \ref{tileT} as 
$$
\mathcal T(S)=  \textup{Thick}_{\varepsilon}\mathcal T(S)\sqcup \bigsqcup_{\sigma\in\mathcal C(S)}\  
\textup{Thin}_{\varepsilon}(\mathcal T(S),\sigma).
$$
Further note that by Proposition \ref{corners} and Proposition \ref{tileT}
$$
\textup{Thin}_{\varepsilon}(\mathcal T(S),\sigma)
\cong\textup{Thick}_{\varepsilon}\mathcal T(S_{\sigma})\times 
\mathbb R^{|\sigma|}\times \mathbb R_{> 0}^{|\sigma|}.
$$

The $\Gamma$-invariant tiling of Teichm\"uller space obtained in Proposition \ref{tileT} induces a
corresponding  tiling of  moduli space as stated in Theorem A in Section 1.

\subsection{A tiling of moduli space: The proof of Theorem A}

 We first show that the canonical projection
$\pi : \mathcal T(S)\rightarrow \mathcal M(S)$ restricted to  any outer cone is injective.
Thus let $CO(x)$ be an outer cone. By its definition it determines a simplex $\sigma\in\mathcal C(S)$ and we can choose adapted FN-coordinates $\Phi_{\tau}$ with $\sigma\subseteq\tau$. For two points $p$ and $q$ in  $CO(x)$ consider their adapted coordinates,
say  $\Phi_{\tau}(p)=(s_*,\theta_*,(l_{\alpha}(p))_{\alpha\in\sigma})$
and $\Phi_{\tau}(q)=(s_*,\theta_*,(l_{\alpha}(q))_{\alpha\in\sigma})$.
If $\pi(p)=\pi(q)$, i.e., $q=h\cdot  p$ for some $h\in \Gamma$, then by Lemma \ref{length}
$$
l_{h^{-1}\cdot \alpha}(p)=l_{\alpha}(h\cdot p)=l_{\alpha}(q)<\varepsilon\ \ \ \textup{for all}\ \ \alpha\in \sigma.
$$
Therefore $l_{\beta}(p)<\varepsilon$ for all $\beta\in h^{-1}\cdot \sigma\cup \sigma$, so that 
$h^{-1}\cdot \sigma\vee \sigma$ is a simplex in $\mathcal C(S)$ (Lemma \ref{disjoint}). We now argue as
in the proof of  Proposition \ref{tileT} to conclude that $h^{-1}\cdot \sigma=\sigma$. Since
$\Gamma$ consists of pure mapping classes we actually have $h^{-1}\cdot \alpha=\alpha$ for  all  $\alpha\in \sigma$. Hence, by the above equalities, $l_{\alpha}( p)=l_{\alpha}(q)$ for  all  $\alpha\in \sigma$,
which then implies that $\Phi_{\tau}(p)=\Phi_{\tau}(q)$ and thus $p=q$. 

Together with the $\Gamma$-invariance of the set of outer cones (see Lemma  \ref{outercone}) this
allows us to define {\it outer cones in moduli space}\,: for $x\in \textup{Thick}_{\varepsilon}\mathcal M(S)$
we set $CO(x):= \pi(CO(\hat{x}))$ where $\hat{x}$ is any lift of $x$.

From Proposition \ref{boundary} and Proposition \ref{tileT} we now get 
$$
\mathcal M(S)=\textup{Thick}_{\varepsilon}\mathcal M(S)\sqcup \bigsqcup_{x\in\partial\textup{Thick}_{\varepsilon}\mathcal M(S)}CO(x)
= \textup{Thick}_{\varepsilon}\mathcal M(S)\sqcup \bigsqcup_{\sigma\in \mathcal E}\  
\bigsqcup_{x\in\mathcal B_{\varepsilon}(\sigma)}CO(x).
$$
Finally, we set
$\textup{Thin}_{\varepsilon}(\mathcal M(S),\sigma):=\{CO(x)\mid x\in\mathcal B_{\varepsilon}(\sigma)\}$.
This completes
 the proof of Theorem A.

\vspace{2ex}

We next reformulate Theorem A in a way which emphasizes the remarkable parallels
with precise reduction theory  for arithmetic groups (see e.g. \cite{OW}, Theorem 3.4 or 
\cite{Sa}, Theorem 9.6). 

Recall that a rational parabolic subgroup of an algebraic group $\mathbf G$ defined over $\mathbb Q$ fixes pointwise a simplex of the rational Tits building of $\mathbf G$ (see e.g. \cite{L1}, Lemma 1.2). Analogously, we call a subgroup $P$ of $\Gamma\subset \textup{Mod}_S$  {\it parabolic} if it fixes pointwise a simplex of the complex of curves $\mathcal C(S)$.
Note that this is in agreement with the characterization of reducible (torsionfree) mapping classes in \cite{Har2}

\begin{corollary}\label{para}  Let $ P_1,\ldots,  P_n$ be representatives of conjugacy classes of the parabolic subgroups   
of  $\Gamma$ (corresponding to simplices $\sigma_1, \ldots,\sigma_n\in \mathcal C(S))$. 
Then there exist subsets $\omega_i$ and $\Omega_i$,  $i=1,\ldots,n$,  of Teichm\"uller space $\mathcal T(S)$, where $\omega_i$ is bounded and $\Omega_i$ is diffeomorphic to $\omega_i\times\mathbb R_{>0}^{|\sigma_i|}$, such that 

\textup{(1)}  The canonical projection $\pi:\mathcal T(S)\longrightarrow \mathcal M(S)$ maps each $\Omega_i$  {\it injectively} into moduli space $\mathcal M(S)$.

\textup{(2)} Each image $\pi(\omega_i)$ in  $\mathcal M(S)$  is compact.

\textup{(3)} For $\varepsilon$ sufficiently small, the moduli space  can be decomposed into a disjoint union 
$$
\mathcal M(S)=\textup{Thick}_{\varepsilon}\mathcal M(S) \sqcup\bigsqcup_{i=1}^n \pi( \Omega_i).
$$ 
\end{corollary}

\emph{Proof}. In view of Theorem A we only have to show the existence of the sets $\omega_i$
and $\Omega_i$. Choose $\varepsilon\leq \varepsilon_0$  as in Proposition \ref{corners}. By Proposition \ref{boundary}  the boundary face $\mathcal B_{\varepsilon}(\sigma_i)$ of 
$\textup{Thick}_{\varepsilon}\mathcal M(S)$ is a torus bundle over $\textup{Thick}_{\varepsilon}\mathcal M(S_{\sigma_i})$ and in particular compact. We can thus choose $\omega_i\subset 
\bigcap_{\alpha\in \sigma_i}\mathcal H_{\alpha}$ as a relatively compact fundamental domain for
the action of  $(\textup{Mod}_{S_{\sigma}}\cap\Gamma)\times \textup{Tw}(\sigma)$ (compare the proof of 
 Proposition \ref{boundary}). Finally we set $\Omega_i:= \bigsqcup_{x\in\omega_i}CO(x)$.
The assertions then follow from Theorem A.
\hfill$\Box$

\vspace{1ex}

Corollary \ref{para} should be compared with the formulation of precise reduction theory
for arithmetic groups given in \cite{BJ}, Proposition III.2.21.
Notice that in this comparison  outer cones are the analoga of Weyl chambers.

\section{The asymptotic cone  of the  moduli space }

In this section we prove Theorem B of Section 1. In particular, we show
that the canonical projection map $\pi : \mathcal T(S)\rightarrow
\mathcal M(S)$ restricted to
certain (coarse) Euclidean cones in $\mathcal T(S)$ is bi-Lipschitz.

The asymptotic cone of the moduli space $\mathcal M(S)$ (see Section 1) is eventually  
identified with a Euclidean cone which is constructed as follows:
We endow each simplex of the complex of curves $\mathcal C(S)$ with a spherical metric and thus 
 obtain
a geometric realization  $|\mathcal C(S)|$  of $\mathcal C(S)$
as a spherical complex.  This spherical metric
 induces a distance function $d_{\mathcal C(S)}$ on the
finite simplical complex $\Gamma\backslash |\mathcal C(S)|$.
The {\it Euclidean cone}  $\textup{Cone}(\Gamma\backslash |\mathcal C(S)|)$  {\it over }
 $\Gamma\backslash |\mathcal C(S)|$ is defined as the product $[0,\infty)\times  
 \Gamma\backslash |\mathcal C(S)|$
with $\{0\}\times \Gamma\backslash |\mathcal C(S)|$ collapsed to a point $\mathcal O$  and  endowed
with the cone metric
$$
d_C^2((a,x),(b,y)) := a^2 + b^2 - 2ab\cos (\min\{\pi, d_{\mathcal C(S)}(x,y)\})  \ \ \textup{(see \cite{BH})} .
$$

Alternatively one can construct $\textup{Cone}(\Gamma\backslash |\mathcal C(S)|)$  as follows.
 Let $\mathcal E$ be the set of simplices of the finite complex $\Gamma\backslash |\mathcal C(S)|$.
 Let $\tau_i\in \mathcal E$, $i=1,\ldots m,$ be the simplices of {\it maximal} dimension $d(S)-1$.
Let $\mathcal B_{\varepsilon_0}(\tau_i)$,  $i=1,\ldots m,$ be the corresponding {\it minimal} boundary faces  of $\partial\textup{Thick}_{\varepsilon_0}\mathcal M(S)$ (which are $d(S)$-dimensional tori, see Proposition \ref{boundary}). For each $i\in\{1,\ldots,m\}$ choose a point $x_i\in \mathcal B_{\varepsilon_0}(\tau_i)$ and consider the closure of the (maximal) outer cone $CO(x_i)\subset \mathcal M(S)$. Each $\overline{CO(x_i)}$ is diffeomorphic to the cone $\mathbb R_{\geq 0}^{d(S)}$ (and actually, by Proposition \ref{metric} below, 
 quasi-isometric to this Euclidean hyperoctant with respect to the 
Teichm\"uller metric). 
The faces of the cones $\overline{CO(x_i)}$ correspond to subsimplices of $\tau_i$. Two cones
$\overline{CO(x_i)}$ and $\overline{CO(x_j)}$ are pasted together along  faces corresponding to subsimplices $\sigma_i\subset
\tau_i$ and $\sigma_j\subset 
\tau_j$ if and only if there is  $h\in \Gamma$ such that $h\cdot 
\sigma_i=\sigma_j$ (the action of $\Gamma\subset \textup{Mod}_S$ being that on the complex of curves).

\subsection{Approximations of the Teichm\"uller metric on thin parts}

We want to study the coarse geometry of moduli space $\mathcal M(S)$ and Teichm\"uller space $\mathcal T(S)$). More
precisely, we are interested in properties of the bi-Lipschitz class of the Teichm\"uller Finsler metric. 
In what follows we describe explicit representatives of that bi-Lipschitz class
on $(\varepsilon,\sigma)$-thin parts of  $\mathcal T(S)$ (and $\mathcal M(S)$). Our approach has
 two key ingredients. The first is a  result of  McMullen. 
He  showed that  $\mathcal M(S)$  
is K\"ahler hyperbolic in the sense of Gromov and thus in particular carries a complete
 finite volume 
 Riemannian metric of
  bounded sectional curvature. The McMullen metric is a modification of  the (incomplete) 
Weil-Petersson metric 
and quasi-isometric to the Teichm\"uller metric (see \cite{Mc}).  
The second key ingredient is  an  expansion of the Weil-Petersson (WP) metric due to  Wolpert 
(see \cite{Wo1}, \cite{Wo2}). 

For each  length function $l_{\alpha}$ we set  $u_{\alpha}:=-\log l_{\alpha}^{1/2}$. 
Considering this logarithmic root
 length instead of $l_{\alpha}$ itself is suggested by  work of  Wolpert (see e.g. 
\cite{Wo1}, \cite{Wo2}). 
Following Wolpert we also set $\lambda_{\alpha}:= \textup{grad}\  l_{\alpha}^{1/2}$
 (resp. 
$\nu_{\alpha}:= -\textup{grad}\  u_{\alpha}$) and define the  \emph{Fenchel-Nielsen-gauge} 
as the differential 1-form 
$\rho_{\alpha}:=
2\pi(l_{\alpha}^{3/2}\langle\lambda_{\alpha},\lambda_{\alpha}\rangle)^{-1}\langle\ ,
J\lambda_{\alpha}\rangle$
$=2\pi(\langle\nu_{\alpha},\nu_{\alpha}\rangle)^{-1}\langle\ ,J\nu_{\alpha}\rangle$ 
(see \cite{Wo2}, 4.15). Note
 that this gauge is
 normalized such that $l_{\alpha}(T_{\alpha})=1$ for $T_{\alpha}:=
(2\pi)^{-1}l_{\alpha}^{3/2}J\lambda_{\alpha}$  the WP-unit infinitesimal FN angle variation. We also set
$\tilde{\rho}_{\alpha}:=l_{\alpha}^{-1/2}\rho_{\alpha}$.

\begin{proposition} \label{metric}There is $\varepsilon_*>0$ depending  only on $S$,
 such that for $\sigma\in\mathcal C(S)$, 
 $\varepsilon \leq\varepsilon_*$  and the   $(\varepsilon,\sigma)$-thin part 
$$
\textup{Thin}_{\varepsilon}(\mathcal T(S),\sigma):=
\bigsqcup_{x\in\bigcap_{\alpha\in\sigma}\mathcal H_{\varepsilon}(\alpha)} CO(x) 
\cong\textup{Thick}_{\varepsilon}\mathcal T(S_{\sigma})\times 
\mathbb R^{|\sigma|}\times \mathbb R_{> 0}^{|\sigma|}
$$ 
the  following Finsler metric expansion of the Teichm\"uller metric with respect to
 adapted FN-coordinates
 holds 
$$
\|.\|^2_{\mathcal T(S)}\asymp 
\|.\|^2_{\mathcal T(S_{\sigma})}
+ \sum_{\alpha\in \sigma}e^{-6u_{\alpha}}\tilde{\rho}_{\alpha}^2+du^2_{\alpha}.
$$
The bi-Lipschitz constants involved in this estimate only depend on $\sigma$ and $\varepsilon_*$.
\end{proposition}

On $(\varepsilon,\sigma)$-thin parts the Teichm\"uller metric is thus quasi-isometric to 
a product metric of a lower dimensional Teichm\"uller space times a product of hyperbolic horoballs.
We remark that Minsky also proved a comparison result for such regions but using  the sup-metric on the products of horoballs  (see \cite{Mi}).

\vspace{1ex}

{\it Proof}. 
McMullen's  modification of the WP metric on Teichm\"uller resp. moduli space is comparable
to the Teichm\"uller metric. More precisely, there is $\varepsilon_1$ (sufficiently small and depending only on $S$) such that for $\varepsilon\leq \varepsilon_1$  and a tangent vector $v\in T_x\mathcal T(S)$ 
the McMullen metric approximation  is given by the formula 
$$
\|v\|^2_{\mathcal T}\asymp \|v\|^2_{WP}+\sum_{l_{\gamma}(x)<\varepsilon}|(\frac{\partial  l_{\gamma}}{l_{\gamma}})(v)|^2,
$$
where the sum is taken over all $\gamma\in\mathcal C(S)$ such that $l_{\gamma}(x)<\varepsilon$
 (see \cite{Mc}, Theorem 1.7). Here the Lipschitz constants only depend on $\varepsilon_1$.
Note that on the $\varepsilon$-thick part there is no modification  and in particular $\|v\|_{\mathcal T}\asymp \|v\|_{WP}$ for $v\in T_x\textup{Thick}_{\varepsilon}\mathcal T(S)$.

 In \cite{Wo2} \S\, 1 and \S\, 4, Wolpert provides an expansion of the Weil-Petersson metric on an $(\varepsilon_2,\sigma)$-thin part $\textup{Thin}_{\varepsilon_2}(S,\sigma)$ considered as a neighbourhood of a point of the positive codimension stratum
 $\mathcal S(\sigma)$ of the Weil-Petersson-completion $\overline{\mathcal T(S)}$ of the Teichm\"uller space. This \emph{augmented
Teichm\"uller space} consists of marked noded Riemann surfaces and is a  stratified (non locally compact) space (see e.g. \cite{Wo1}). The strata are indexed by the simplices of $\mathcal C(S)$ and are (products of) lower dimensional Teichm\"uller spaces. A boundary stratum $\mathcal S(\sigma)\subset \overline{\mathcal T(S)}$ consists
of those marked Riemann surfaces whose nodes correspond bijectively to the vertices of $\sigma$
and is isomorphic 
to the Teichm\"uller space $\mathcal T(S_{\sigma})$ of the cut surface $S_{\sigma}$ (obtained 
by cutting $S$ along the circles of $\sigma$). Using the FN-gauges $\rho_{\alpha}$ Wolpert's expansion of the WP-metric 
can be written as
\begin{eqnarray*}
 \langle\ ,\ \rangle_{WP} & \asymp & \langle\ ,\ \rangle_{WP}^{\mathcal T(S_{\sigma})}
+\sum_{\alpha\in\sigma}|\partial l_{\alpha}^{1/2}|^2\\
 & \asymp & \langle\ ,\ \rangle_{WP}^{\mathcal T(S_{\sigma})}
+\sum_{\alpha\in\sigma}(dl^{1/2}_{\alpha}\circ J)^2+(dl^{1/2}_{\alpha})^2
\asymp \langle\ ,\ \rangle_{WP}^{\mathcal T(S_{\sigma})}
+\sum_{\alpha\in\sigma}l_{\alpha}^3\rho^2_{\alpha}+(dl^{1/2}_{\alpha})^2.
\end{eqnarray*}
The involved 
Lipschitz constants depend on $\sigma$ and  $\varepsilon_2$.

We now apply these estimates to the $(\varepsilon,\sigma)$-thin part 
$$\textup{Thin}_{\varepsilon}(\mathcal T(S),\sigma)
\cong\textup{Thick}_{\varepsilon}\mathcal T(S_{\sigma})\times 
\mathbb R^{|\sigma|}\times \mathbb R_{> 0}^{|\sigma|}.$$
 Recall that $\|.\|_{\mathcal T}\asymp \|.\|_{WP}$ on $\textup{Thick}_{\varepsilon}\mathcal T(S_{\sigma})$.
We set $\varepsilon_*:=\min(\varepsilon_1,\varepsilon_2)$ and let $\varepsilon\leq\varepsilon_*$. Inserting
the McMullen formula in the above expansion and using that 
  $|\partial l_{\gamma}(v)|^2= (dl_{\gamma}(v))^2+(dl_{\gamma}(Jv))^2$ and
$dl_{\gamma}=2 l^{1/2}_{\gamma} dl^{1/2}_{\gamma}$ we obtain 
\begin{eqnarray*}
\|.\|^2_{\mathcal T(S)} & \asymp &  \|.\|^2_{\mathcal T(S_{\sigma})}+\sum_{\alpha\in\sigma} (dl_{\alpha}^{1/2}\circ J)^2(1+\frac{1}{l_{\alpha}})+
(dl_{\alpha}^{1/2})^2 (1+\frac{1}{l_{\alpha}})\\
   & \asymp & \|.\|^2_{\mathcal T(S_{\sigma})}+\sum_{\alpha\in\sigma}
\frac{1}{l_\alpha}(dl^{1/2}_{\alpha}\circ J)^2+
\frac{1}{l_\alpha}(dl^{1/2}_{\alpha})^2\ \ \ \ (\textup{since}\  l_{\alpha}\leq \varepsilon)\\
& \asymp & \|.\|^2_{\mathcal T(S_{\sigma})}+\sum_{\alpha\in\sigma} l_{\alpha}^2\rho_{\alpha}^2+
\frac{1}{l_\alpha}(dl^{1/2}_{\alpha})^2\ \ \ \  (\textup{by Wolpert's expansion above).}
\end{eqnarray*}
We next set 
$\tilde{T}_{\alpha}:=l_{\alpha}^{1/2}T_{\alpha}$ for $T_{\alpha}$ the WP-unit infinitesimal FN-angle variation.
Then $\|\tilde{T}_{\alpha}\|_{\mathcal T(S)}\asymp l_{\alpha}^{-1/2}\|\tilde{T}_{\alpha}\|_{WP}=1$ and hence
$\tilde{\rho}_{\alpha}=l_{\alpha}^{-1/2}\rho_{\alpha}$ is the renormalized FN gauge: $\tilde{\rho}_{\alpha}(\tilde{T}_{\alpha})=
\rho_{\alpha}(T_{\alpha})=1$.
Substituting  this together with  $u_{\alpha}=-\log l^{1/2}_{\alpha} $ in the previous expression 
we eventually get the claimed estimate
$$
\|.\|^2_{\mathcal T(S)}  \asymp
\|.\|^2_{\mathcal T(S_{\sigma})}+\sum_{\alpha\in\sigma}e^{-6u_{\alpha}}\tilde{\rho}^2_{\alpha}+
 du_{\alpha}^2.
 $$
The Lipschitz constants of this comparison depend only on $\sigma$ and $\varepsilon_*$.
\hfill $\Box$

\vspace{1ex}

\begin{corollary} \label{finslercover} 
The 
 Finsler metric approximation on the $(\varepsilon,\sigma)$-thin parts 
$\textup{Thin}_{\varepsilon}(\mathcal T(S),\sigma)$ of $\mathcal T(S)$ given in 
Proposition \ref{metric} descends
to the $(\varepsilon,\sigma)$-thin parts $\textup{Thin}_{\varepsilon}(\mathcal M(S),\sigma)\cong \mathcal B_{\varepsilon}(\sigma)\times \mathbb R_{>0}^{|\sigma|}$
of $\mathcal M(S)$.
\end{corollary}

\vspace{1ex}

{\it Proof}.\ 
The canonical projection $\pi:\mathcal T(S)\longrightarrow\mathcal M(S)$
is a Finsler covering for the (induced) Teichm\"uller Finsler metrics. By Proposition \ref{boundary}
and Theorem A (and the subresults of their proofs) we further have
$$
\textup{Thin}_{\varepsilon}(\mathcal M(S),\sigma)=(\textup{Mod}_{S_{\sigma}}\cap\Gamma)\times \textup{Tw}(\sigma)\backslash\textup{Thin}_{\varepsilon}(\mathcal T(S),\sigma).
$$
Finally, note that the group $(\textup{Mod}_{S_{\sigma}}\cap\Gamma)\times\textup{Tw}(\sigma)$ leaves the factorization 
$$\textup{Thin}_{\varepsilon}(\mathcal T(S),\sigma)
\cong \textup{Thick}_{\varepsilon}\mathcal T(S_{\sigma})\times\mathbb R^{|\sigma|}\times \mathbb R_{> 0}^{|\sigma|} $$
used in Proposition \ref{metric} invariant. \hfill $\Box$

\begin{corollary} \label{cmetric} For adapted FN-coordinates $(s,\theta,u)$ and $\varepsilon > 0$ as in Proposition  \ref{metric} define
the projection map
$$
\pi_u: \textup{Thin}_{\varepsilon}(\mathcal T(S),\sigma)\cong \textup{Thick}_{\varepsilon}T(S_{\sigma})\times \mathbb R^{|\sigma|}\times \mathbb R_{> 0}^{|\sigma|}\longrightarrow \mathbb R_{> 0}^{|\sigma|};\ \ (s,\theta,u)\longmapsto u.
$$
Then $\pi_u$ is length non-increasing on the the  $(\varepsilon,\sigma)$-thin part $\textup{Thin}_{\varepsilon}(S,\sigma) $
 with respect to the restrictions of the approximate Teichm\"uller metric.
On $\mathbb R_{> 0}^{|\sigma|}$ this restriction is comparable with the Euclidean metric. 
\end{corollary}

\subsection{Bi-Lipschitz images of  Euclidean cones and a net in moduli space}

Throughout this section we
choose a fixed positive real number $\varepsilon_0$ small enough such that
both Proposition \ref{corners}   and Proposition \ref{metric}  hold.
As usual we denote by
$\pi : \mathcal T(S)\rightarrow \mathcal M(S)$   the canonical projection.
The metric on $\mathcal M(S)$ is given  by
$d_{\mathcal M}(\pi(x), \pi(y)) = \inf_{h\in\Gamma} d_{\mathcal T}(x, h\cdot y)$.
Recall from Section 2.3 that outer cones in $\mathcal M(S)$ are projections of
outer cones in $\mathcal T(S)$.

We next point out metric properties of the  tiling in Theorem A.

\begin{proposition}\label{bilip}
Let $\varepsilon\leq\varepsilon_0$ be as in Proposition \ref{corners} and let $CO(x)$,
$x\in \mathcal B_{\varepsilon}(\sigma) $, be an outer cone of $\textup{Thick}_{\varepsilon}\mathcal M(S)$.
Then  the restriction of the canonical projection $\pi : \mathcal T(S)\rightarrow \mathcal M(S)$
to any lifted outer cone $CO(z)\subset \mathcal T$, with $\pi(z)=x$, is bi-Lipschitz with respect
 to the Teichm\"uller distance functions of $\mathcal T(S)$ and $\mathcal M(S)$, respectively.
\end{proposition}

{\it Proof}.  In a lift
$CO(z)$ of $CO(x)$ to $\mathcal T(S)$  we pick two arbitrary points $v,w$.
 We wish to show that $d_{\mathcal M}(\pi(v),\pi(w))\succ d_{\mathcal T}(v,w)$
(the opposite inequality holds
by definition). Let $c$ be a  curve of minimal length in $\mathcal M(S)$ between $\pi(v) $ and $\pi(w)$, i.e.,
 $L(c)=d_{\mathcal M(S)}(\pi(v),\pi(w))$.
We distinguish two cases: (1) $c$ does {\it not} intersect the compact $\varepsilon$-thick part $\textup{Thick}_{\varepsilon}\mathcal M(S)$ and 
(2) $c$ intersects 
$\textup{Thick}_{\varepsilon}\mathcal M(S)$.

Case (1):  According to Theorem  A  we can decompose the curve
$$
c: [0,d_{\mathcal M(S)}(\pi(v),\pi(w))]\longrightarrow \mathcal M(S)\setminus \textup{Thick}_{\varepsilon}\mathcal M(S)
$$
 into a finite
number of segments $c_i, 1\leq i\leq l$, say, each of which is contained in a subset
of $\mathcal M(S)$ diffeomorphic to  $\mathcal B_{\varepsilon}(\sigma)\times \mathbb R_{> 0}^{|\sigma|}$
where  $\mathcal B_{\varepsilon}(\sigma)$ is a torus 
bundle over the $\varepsilon$-thick part of a moduli space and compact.  
 By its definition each  outer cone of $\mathcal B_{\varepsilon}(\sigma)$ is diffeomorphic  to the open
 hyperoctant $\mathbb R_{>0 }^{|\sigma|}$. 

By  Proposition \ref{metric} and Corollary \ref{finslercover}
a representative of  the bi-Lipschitz class of the Teichm\"uller  metric
 at  points
$(b,u)\in \mathcal B_{\varepsilon}(\sigma)\times \mathbb R_{> 0}^{|\sigma|}$ can be written in the 
form 
$$
ds_{(b,u)}^2\asymp db_{(b,u)}^2+du_u^2,
$$
where $du_u^2$ is the Euclidean(!) metric  on $ \mathbb R_{> 0}^{|\sigma|}$. 
For the length of the segment $c_i(s)=(b_i(s),u_i(s)), s\in [s_i,s_{i+1}]$, we thus have
the estimate $L(c_i)\succ L(u_i), 1\leq i\leq l$ (compare also Corollary \ref{cmetric}). Hence
$$
d_{\mathcal M(S)}(\pi(v),\pi(w))=L(c)=\sum_{i=1}^l L(c_i)\succ  \sum_{i=1}^l L(u_i).
$$

In view of the second construction of $\textup{Cone}(\Gamma\backslash |\mathcal C(S)|)$  each outer cone of $\textup{Thick}_{\varepsilon}\mathcal M(S)$
can  be mapped diffeomorphically to a Euclidean cone in  $ \textup{Cone}(\Gamma\backslash |\mathcal C(S)|)$. 
Denote by $\overline{v}$ and 
$\overline{w}$ the corresponding images of $v,w$ in  $ \textup{Cone}(\Gamma\backslash |\mathcal C(S)|)$. 
Then the continuous curve $u:=u_{1}\ast u_2\ast \cdots \ast u_l$ can be indentified with a curve in $ \textup{Cone}(\Gamma\backslash |\mathcal C(S)|)$, which 
connects $\overline{v}$ to
$\overline{w}$   and thus satisfies
$$
\sum_{i=1}^l L(u_i)\geq d_C(\overline{v},\overline{w}).
$$

By Theorem A  we 
have $CO(z)\cong CO(x)$ and $d_{C}(\overline{v}, \overline{w})=d_{Euc}(\pi(v), \pi(w))=d_{Euc}(v, w)$. 
Finally, by Proposition \ref{metric}, $d_{Euc}(v,w)\succ d_{\mathcal T}(v,w)$.

Together with the previous estimate this proves the claim in case (1).

Case (2): We decompose $c$ into three segments $c=c_1\ast c_2\ast c_3$, such that $c_1$ connects
$\pi(v)$ to the first
intersection point of $c$ with $\textup{Thick}_{\varepsilon}\mathcal M(S)$ and  $c_3$ connects the last intersection point
of $c$ with $\textup{Thick}_{\varepsilon}\mathcal M(S)$ to $\pi(w)$. Since $c$ is minimal we then have
$$
d_{\mathcal M}(\pi(v),\pi(w))=L(c)=L(c_1)+L(c_2)+L(c_3)\geq L(c_1)+L(c_3).
$$
We can now argue as in case (1). We consider the $u$-projections $u_1$ of $c_1$ (resp. $u_3$ of $c_3$) and their
images $\overline{u}_1$ (resp. $\overline{u}_3$) in $\textup{Cone}(\Gamma\backslash |\mathcal C(S)|)$. Let also $\overline{v}$ 
(resp. $\overline{w}$) be the images of $\pi(v)$ (resp. $\pi(w)$) in $\textup{Cone}(\Gamma\backslash |\mathcal C(S)|)$.
Then $\overline{u}_1$ (resp. $\overline{u}_3$) joins $\overline{v}$ to the apex  $\mathcal O$ (resp. $\mathcal O$ to $\overline{w}$).
 Hence as in case (1)

\begin{eqnarray*}
L(c_1)+L(c_3) & \succ & L(u_1)+L(u_3)\geq  d_C(\overline{v},\mathcal O) + d_C(\mathcal O,\overline{w})\geq d_C(\overline{v},\overline{w})=\\
 & = & d_{Euc}(\pi(v),\pi(w))
= d_{Euc}(v,w)\succ d_{\mathcal T}(v,w).
\end{eqnarray*}
This proves the claim $d_{\mathcal M}(\pi(v),\pi(w)) \succ d_{\mathcal T}(v,w)$ also in case (2).
\hfill $ \Box$
\vspace{2ex}

 Recall that a
 subset $\mathcal N$ of a metric
space $(X,d)$ is called a ($\delta$--){\it net} if there is a positive  constant
$\delta$ such that  $d(p,\mathcal N)\leq \delta$ for all $p\in X$; in particular the
Hausdorff-distance
between $\mathcal N$ and $X$ is at most $\delta$.

\begin{proposition} \label{net} There is a net $\mathcal N$  in moduli space $\mathcal M(S)$ consisting of finitely  many
 bi-Lipschitz embedded  $d(S)$-dimensional Euclidean cones. In particular, if $d(S)\geq 2$, then
$\mathcal M(S)$ is not Gromov-hyperbolic.\end{proposition}

{\it Proof}. Let $\varepsilon$ be as in Proposition \ref{bilip}. By Theorem A and since 
$\textup{Thick}_{\varepsilon}\mathcal M(S)$ is compact and thus has finite diameter,
$$\bigsqcup_{\sigma\in \mathcal E}\textup{Thin}_{\varepsilon}(\mathcal M(S),\sigma)\cong
\mathcal B_{\varepsilon}(\sigma)\times \mathbb R_{>0}^ {|\sigma|}
$$ 
is a net in $\mathcal M(S)$.
Note that $\mathcal B_{\varepsilon}(\sigma)$ is compact and that on $\textup{Thin}_{\varepsilon}(\mathcal M(S),\sigma)$ one has $u_{\alpha}\geq -\log \varepsilon^{1/2}$, $\alpha\in\sigma$.
Hence  the metric expansion in Corollary \ref{finslercover}
implies that 
$$
d_{\mathcal M}((b_1,u), (b_2,u))\prec d_{\mathcal B}(b_1,b_2)\leq \textup{const}(\varepsilon,\sigma) \ \ (b_1,b_2\in\mathcal B_{\varepsilon}(\sigma); u\in \mathbb R_{>0}^ {|\sigma|}).
$$
Therefore any cone $CO(x_{\sigma})\subset\textup{Thin}_{\varepsilon}(\mathcal M(S),\sigma), x_{\sigma}\in
\mathcal B_{\varepsilon}(\sigma)$, is a net
in $\textup{Thin}_{\varepsilon}(\mathcal M(S),\sigma)$. Recall that $\mathcal E$ denotes the finite 
set of simplices in $\Gamma\backslash |\mathcal C(S)|$.
Given some $\sigma\in \mathcal E$ there is $\tau\in\mathcal E$ of maximal dimension $d(S)-1$ containing $\sigma$ and by Lemma \ref{bdpoint-simplex}, $\mathcal B_{\varepsilon}(\tau)\subset \mathcal B_{\varepsilon}(\sigma)$.
We can thus choose the above cone $CO(x_{\sigma})\cong \mathbb R_{>0}^{|\sigma|}$ 
as a face of a closed cone $\overline{CO(x_{\tau})}$ for some $x_{\tau}\in \mathcal B_{\varepsilon}(\tau)$. The net 
 $\mathcal N$ can now be constructed as follows. 
 Let $\tau_i\in \mathcal E$, $i=1,\ldots m,$ be the simplices of {\it maximal} dimension $d(S)-1$.
Let $\mathcal B_{\varepsilon}(\tau_i)$,  $i=1,\ldots m,$ be the corresponding {\it minimal} boundary faces  of $\partial\textup{Thick}_{\varepsilon}\mathcal M(S)$ (which are $d(S)$-dimensional tori, see Proposition \ref{boundary}). For each $i\in\{1,\ldots,m\}$ we choose a point $x_i\in \mathcal B_{\varepsilon}(\tau_i)$ and  then we set  $\mathcal N:=\bigcup_{i=1}^m\overline{CO(x_i)}$. \hfill $\Box$

\subsection{A quasi-isometric approximation of moduli space}

We recall the notion of Gromov-Hausdorff-convergence of (unbounded) pointed metric spaces
 (see \cite{Gr2}, Chapter 3). First,  the {\it distortion} of a map $f: A\rightarrow B$ of metric spaces $A$ and $B$  is defined
as
$$\textup{dis}(f) := \sup_{a,b}|d_A(a,b) - d_B(f(a), f(b))|.
$$
 The {\it uniform distance}
between metric spaces $A$ and $B$ is defined as $\textup{u-dist}(A,B) := \inf_f\textup{dis}(f)$ where the
infimum is taken over all bijections $f:A\rightarrow B$.
A sequence of metric
 spaces $X_n$ {\it Hausdorff-converges} to a metric space $X$ if and only if for
every $\delta>0$ there is an $\delta$--net $X_{\delta}$ in $X$ which is the
uniform limit of $\delta$--nets $(X_n)_{\delta}$ in $X_n$.
We say that a sequence $(X_n,p_n)$ of unbounded, pointed metric spaces Hausdorff-converges
to a pointed
metric space $(X,p)$ if for every $r>0 $ the balls $B_r(p_n)$ in $X_n$  Hausdorff-converge
to the ball $B_r(p)$ in $X$. 

Let $x_0$ be an (arbitrary)  point of the moduli space $\mathcal M(S)$.
The {\it asymptotic cone} of $\mathcal M(S)$ endowed with the Teichm\"uller distance $d_{\mathcal M}$
is defined as the Gromov-Hausdorff-limit of pointed metric spaces:
$$ \textup{As-Cone}(\mathcal M(S)) := {\mathcal H}-{\lim}_{n\rightarrow\infty} (\mathcal M(S), x_0, \frac{1}{n} d_{\mathcal M}).$$

By  Proposition \ref{net} there is a net $\mathcal N$ \ in $\mathcal M(S)$ consisting of disjoint closed quasi-Euclidean (outer) cones
$\overline{CO(x_i)}$, $i=1,\ldots,m$.  
FN-coordinates $(u_{\alpha}=-\log l_{\alpha}^{1/2})_{\alpha\in\tau_i}$  give rise to  diffeomorphisms $\Phi_i:\overline{CO(x_i)}\cong \mathbb R_{\geq a}^{d(S)}\cong \mathbb R_{\geq 0}^{d(S)}$, $i=1,\ldots,m$ and $a=-\log\varepsilon^{1/2}$.  Using the second construction of $\textup{Cone}(\Gamma\backslash |\mathcal C(S)|)$ we have a map
$$
f_1: \mathcal N \rightarrow \textup{Cone}(\Gamma\backslash |\mathcal C(S)|);\  \  (u_{\alpha})_{\alpha\in\sigma_i}\longmapsto (\Phi_i(u_{\alpha}))_{\alpha\in\sigma_i}.
$$ 
 For any $n\in \mathbb N$ we then define a map
$$f_n: \mathcal N\subset (\mathcal M(S), \frac{1}{n}d_{\mathcal M})\rightarrow
\textup{Cone}(\Gamma\backslash |\mathcal C(S)|);\ \
f_n(x) := \frac{1}{n} f_1(x).
$$

By Proposition \ref{bilip} (or Theorem A)
 $f_n$ is a
bijection from the interior of $\mathcal N$ onto its image  in $\textup{Cone}(\Gamma\backslash |\mathcal C(S)|)$. This image  is open and dense in $\textup{Cone}(\Gamma\backslash |\mathcal C(S)|)$ for all $n$.

\begin{proposition} \label{qi}There is a constant $D> 0$ such that
$\textup{dis}(f_n)\leq \frac{1}{n}D$ for all $n\in \mathbb N$, i.e., $f_n$ is a
quasi-isometry. In particular,  $\mathcal M(S)$ is quasi-isometric to
$\textup{Cone}(\Gamma\backslash |\mathcal C(S)|)$.
\end{proposition}

{\it Proof}. We consider a polyhedral exhaustion $\mathcal M (S) = \bigcup_{\varepsilon\leq \varepsilon_0}\textup{Thick}_{\varepsilon}\mathcal M(S)$ (see Corollary \ref{exhaustion}). The
intersection $\mathcal N\,'$ of the net $\mathcal N$ with  $\mathcal M(S)\setminus \textup{Thick}_{\varepsilon_0}\mathcal M(S)$ is still a net in $\mathcal M(S)$ since $\textup{Thick}_{\varepsilon_0}\mathcal M(S)$ is
compact. Let $u,v$ be two points  in the interior of $\mathcal N\,'$.
We take a path $c([0,L])$ in $\mathcal M(S)$ between $u$ and $v$ of {\it minimal} length and parametrized by
arc-length. Since the diameter of $\textup{Thick}_{\varepsilon_0}\mathcal M(S)$ is finite we can assume that $c$ lies in $\mathcal M(S)\setminus \textup{Thick}_{\varepsilon_0}\mathcal M(S)$ (taking an additive constant into account).  By Theorem A there is a tiling of $\mathcal M(S)\setminus \textup{Thick}_{\varepsilon_0}\mathcal M(S)$ into finitely many $(\varepsilon_0,\sigma)$-thin parts  $\mathcal M_i\cong\mathcal B_{\varepsilon_0}(\sigma_i)\times \mathbb R_{> 0}^{k_i}, i=1,\ldots,q$. We can thus find $t_i\in[0,L]$, $ 0\leq i\leq n+1$ such
that
$t_0 = 0, t_{n+1}= L$ and   $[0,L] = \bigcup_{i=0}^n[t_i,t_{i+1}]$ with $c([t_i,
t_{i+1}])\subset \mathcal M_{j_i}$ for $j_i\in\{1,\ldots,q\}$. Associated to the path
$c$ we thus get a well-defined string
$(j_0,j_1,\ldots,j_n)$.
We next replace $c$ by a path $\overline{c}$ whose associated string contains an element
$j_k\in  \{1,\ldots,q\}$ at most once: Start with $j_0$ and assume that $0\leq l\leq n$ is
the greatest index
such that $j_l = j_0$. Geometrically this means that the path $c$ returns to
$\mathcal M_{j_0}$ at $c(t_l)$. By Proposition \ref{metric} and Proposition \ref{bilip} the shortest curve between two points 
in $\mathcal M_i=\mathcal B_{\varepsilon_0}(\sigma_i)\times \mathbb R_{> 0}^{k_i}$ is coarsly
equivalent to a straight line
in $\mathbb R_{> 0}^{k_i}$.
We can thus  replace $c([t_0,t_{l+1}])$ by a  segment, say
$\overline{c}([s_{j_0},s_{j_0+1}])$, in $ \mathbb R_{> 0}^{k_i}$ whose length $L$  satisfies
$$L(\overline{c}([s_{j_0},s_{j_0+1}]))\leq L(c[t_0,t_{l+1}]) + 2D_1,$$
for  $D_1:=\max_{1\leq i\leq m}\textup{diam}(\mathcal B_{\varepsilon_0}(\sigma_i))$.  Repeating this procedure  we eventually
get
a sequence of at most $q$ segments
$\overline{c}([s_k,s_{k+1}])$ in $\mathcal N\,'$ of total length $\leq L(c) + 2qD_1.$
Moreover we can choose $\overline{c}([s_k,s_{k+1}])$ in such a way that the endpoint
$\overline{c}(s_{k-1})\in \mathcal M_{q_{k-1}}$ and the initial point
$\overline{c}(s_{k})\in \mathcal M_{q_{k}}$ are on the same levelset of the exhaustion
 of $\mathcal M(S)$ we are considering.
 
We next show that these segments can be modified in such a way that they are mapped by $f_1$ to a
continous path
$\tilde{c}$ in $\textup{Cone}(\Gamma\backslash |\mathcal C(S)|)$ from $f_1(u)$ to  $f_1(v)$.
With respect to the representative of the bi-Lipschitz class of the Teichm\"uller  metric
described in Proposition \ref{metric} each outer cone $CO(x_i)$ in $\mathcal N$ is isometric to a Euclidean
simplicial cone in $\textup{Cone}(\Gamma\backslash |\mathcal C(S)|)$. The map $f_1$ thus yields  well defined
images of (the interior of)  the segments
$\overline{c}([s_k,s_{k+1}])\in \mathcal M_{q_k}$. 
By construction the endpoint
$\overline{c}(s_{k-1})\in \mathcal M_{q_{k-1}}$ and the initial point
$\overline{c}(s_{k})\in \mathcal M_{q_{k}}$ are on the same levelset of the exhaustion
  and are at most the  distance $2D_1$ apart.
By the (second) construction of $\textup{Cone}(\Gamma\backslash |\mathcal C(S)|)$ we can join their  images
 by segments in $\textup{Cone}(\Gamma\backslash |\mathcal C(S)|)$ of uniformly bounded length
 $2D_1$ to obtain the path $\tilde{c}$. This argument yields the
 estimate
$$
d_C(f_1(u),f_1(v)) \leq L_C(\tilde{c})\leq 4qD_1+ L(c)= 4qD_1+ d_{\mathcal M}(u,v).
$$

On the other hand, given a geodesic path  $\tilde{c}$ in $\textup{Cone}(\Gamma\backslash |\mathcal C(S)|)$
joining
two  points $x$ and $y$ in the image under
$f_1$ of the interior of  $\mathcal M(S)\setminus \textup{Thick}_{\varepsilon_0}\mathcal M(S)$, we can lift the segments contained in the simplicial cones
$C(\tau_j)$, say, in $\subset \textup{Cone}(\Gamma\backslash |\mathcal C(S)|)$  to the corresponding outer cones $CO(x_j)$ via $f_1^{-1}$. By
the same
arguments as before the endpoints of the lifted segments can be joined in $\mathcal M(S)$ to form a
continuous path $c$ between $f_1^{-1}(x)$ and $f_1^{-1}(y)$ of length
$$
d_{\mathcal M}(f_1^{-1}(x),f_1^{-1}(y)) \leq L(c) \leq  2qD_1 + L(\tilde{c}) \leq 2qD_1 + d_C(x,y).
$$
Combining $(1)$ and $(2)$  and setting $D:=4qD_1$ we get for all $u,v$ in (the interior of)
$\mathcal M(S)\setminus \textup{Thick}_{\varepsilon_0}\mathcal M(S)$
that
$$|d_{\mathcal M}(u,v) - d_C(f_1(u), f_1(v))|\leq D.$$
Finally,  since for any $n\in\mathbb N$
\begin{eqnarray*}
|\frac{1}{n}d_{\mathcal M}(u,v) - d_C(f_n(u), f_n(v))| & = & |\frac{1}{n}d_{\mathcal M}(u,v) -
\frac{1}{n}d_C(f_1(u), f_1(v))| = \\
 & =&  \frac{1}{n} |d_{\mathcal M}(u,v) - d_C(f_1(u), f_1(v))| \leq \frac{D}{n},
 \end{eqnarray*}
which completes the proof of Proposition \ref{qi}.
\hfill $\Box$

\subsection{The proof of Theorem B}

 Given $\delta >0$ there is $n_{\delta}\in\mathbb N$, such that
$\mathcal N\,' := \mathcal N\cap(\mathcal M(S)\setminus \textup{Thick}_{\varepsilon_0}\mathcal M(S))\subset \mathcal M(S)$ is an $\delta$--net in $(\mathcal M(S), \frac{1}{n}d_{\mathcal M})$
for all $n\geq n_{\delta}$.
Since $f_1(\mathcal N\,')$ is dense in $\textup{Cone}(\Gamma\backslash |\mathcal C(S)|)\setminus f_1(\mathcal N\cap \textup{Thick}_{\varepsilon_0}\mathcal M(S))$ it is a
$\delta$--net in the latter set for all
$\delta>0$. Then  there is $r_0$ such that for  $r\geq r_0$  the same assertions are
true
for the subsets of balls with radius $r$:
$$
\mathcal N\,'\cap B_r(v_0)\subset (\mathcal M(S), \frac{1}{n}d_{\mathcal M})\ \textup{and}\  f_n(\mathcal N\,')\cap B_r(\mathcal O)
\subset (\textup{Cone}(\Gamma\backslash|\mathcal C(S)|), d_C).
$$
From Proposition \ref{qi} and its proof we obtain the following uniform estimate:
$$
\textup{u-dist}(\mathcal N\,'\cap B_r(v_0),f_n(\mathcal N')\cap B_r(\mathcal O))\leq \frac{D}{n},
$$
which completes the proof of Theorem B.

\begin{corollary} The asymptotic cone of moduli space $\mathcal M(S)$ endowed with the Teichm\"uller metric  is bi-Lipschitz equivalent  to a $\textup{CAT(0)}$ space and in particular 
contractible. 
\end{corollary}

{\it Proof}.
 By Theorem B, $\textup{As-Cone}(\mathcal M(S))$ is bi-Lipschitz equivalent to the Euclidean cone
$\textup{Cone}(\Gamma\backslash |\mathcal C(S)|)$. A theorem of  Berestovski says that for a geodesic
metric
space $Y$ the Euclidean cone $C(Y)$ is a CAT(0) space  if and only if $Y$ is a CAT(1) space (see
[9], Ch. II.4). Since $\Gamma\backslash |\mathcal C(S)|$ is obtained 
by glueing spherical simplices it is CAT(1) (see \cite{BH}, Ch II, 4.4).
\hfill$\Box$

\section{Positive scalar curvature}

 Block and  Weinberger showed that  on a locally symmetric space $\Gamma\backslash G/K$
there exists a complete Riemannian metric of uniformly positive scalar curvature  if and only if 
$\Gamma$ is an arithmetic group of $\mathbb Q$-rank $\geq 3$ (see \cite{BW}). These Riemannian
metrics have the quasi-isometry type of  rays and one may ask whether there is also a metric 
that  is both uniformly positively curved and  quasi-isometric to the symmetric metric. That question was 
answered negatively  for locally symmetric spaces associated to arithmetic groups 
of $\mathbb Q$-rank $\geq 2$ by  Chang (see \cite{Ch}).

Farb and Weinberger recently anounced a result  analogous to \cite{BW} for $\mathcal M(S)$: A finite cover of the
moduli space admits a complete, finite volume Riemannian metric of uniformly bounded positive 
scalar curvature if and only if $d(S)\geq 3$ (see \cite{Fa}, 4.2). Again this metric has  the quasi-isometry type of  rays
and thus, by Theorem B,  is not quasi-isometric to the Teichm\"uller metric. Farb conjectured  the analogue of 
Chang's result (see \cite{Fa}, Conjecture 4.6). Our Theorem C in Section 1 confirms this conjecture.

\subsection{Outline of the proof of Theorem C}

Given Theorem A and the subresults of its proof, the proof of Theorem  C is along exactly the same lines as that of Theorem 1 in \cite{Ch}. 
We therefore only sketch those  ingredients which rely on general results and  in particular refer to \cite{Ch} for definitions
and more details.

 Given a tiling of
the moduli space  $\mathcal M(S)$ as in Theorem A we consider a (maximal) outer cone
$CO(x)$ of $\partial \textup{Thick}_{\varepsilon_0}\mathcal M(S) $ corresponding to a  simplex $\tau$ 
of maximal dimension in the curve complex $\mathcal C(S)$.
Then  $CO(x)$ is contained in a subset of $\mathcal M(S)$ diffeomeorphic to $\mathcal B_{\varepsilon_0}(\tau)\times \mathbb R_{\geq 0}^{d(S)}$,
where (in this case) $\mathcal B_{\varepsilon_0}(\tau)\cong \textup{Tw}(\tau)\backslash \mathbb R^{d(S)}\cong T^{d(S)}$ is a $d(S)$-dimensional torus (see Proposition \ref{boundary}). Using adapted FN-coordinates $(\theta,u)$
on $\mathcal B_{\varepsilon_0}(\tau)\times\mathbb R_{\geq 0}^{d(S)}$ we consider
the projection 
$$
\pi_u:\mathcal B_{\varepsilon_0}(\tau)\times\mathbb R_{\geq 0}^{d(S)}\longrightarrow \mathbb R_{\geq 0}^{d(S)};
\ \ \ (\theta,u)\longmapsto u.$$

We then select a hypersurface $Q$ in the hyperoctant $\mathbb R_{\geq 0}^{d(S)}$ sufficiently
far from the faces, so that the inverse image under $\pi_u$ of each point of $Q$ is a $d(S)$-dimensional torus.
Then $V:=\pi_u^{-1}(Q)$ is diffeomorphic to $T^{d(S)}\times\mathbb R^{d(S)-1}$.
The complement of the hypersurface $V$ in $\mathcal M(S)$ consists of two components. Let $A$ be the closure of the  component
containing the fibre $\pi_u^{-1}(0)$ and let $B$ be the closure 
of $\mathcal M(S)\setminus \textup{int}\, A$. Then the pair $(A,B)$ forms a ``coarsely excicive decomposition'' of 
$\mathcal M(S)$  with intersection $A\cap B=V$ (see \cite{Ch} for the definition).

One can now proceed exactly as in the proof of Theorem 1 in \cite{Ch}:
\begin{itemize}

\item First note that
$\pi_1(V)\cong \textup{Tw}(\tau)\cong \mathbb Z^{d(S)}$. Hence there is an injection 
$\pi_1(V)\hookrightarrow \pi_1(\mathcal M(S))=\Gamma$ and a K-theoretic Mayer-Vietoris 
sequence applies. 

\item Next let $D$ be the 
lifted classical Dirac operator on the pullback spinor bundle of $\mathcal T(S)$. Then one defines 
a generalized Roe algebra $C_{\Gamma}^*\mathcal M(S)$ and a generalized coarse index $\textup{Ind}(D)\in K_*(C_{\Gamma}^*\mathcal M(S))$
as in \cite{Ch}.

\item Finally, the fact that the strong Novikov conjecture for nilpotent groups is true,     
 implies that
$\textup{Ind}(D)\neq 0$.  Chang proved in \cite{Ch} that this is an obstruction to  the existence of a uniformly bounded positive scalar
curvature metric in the same quasi-isometry class as the given metric (which in the present case is the Teichm\"uller metric).
\end{itemize}

\vspace{4ex}

\noindent\textsc{Institute for Algebra and Geometry\\
 University of Karlsruhe, 76131 Karlsruhe, Germany}

\vspace{1ex}

\noindent{\tt enrico.leuzinger@math.uni-karlsruhe.de}
\end{document}